\documentclass[12pt,a4paper]{amsart}

\usepackage[T1]{fontenc}
\usepackage{lmodern}
\usepackage{amsmath,amssymb,mathtools}
\usepackage[margin=1in]{geometry}
\usepackage{microtype}
\usepackage{booktabs}
\usepackage{array}
\usepackage{enumitem}
\usepackage{needspace}
\usepackage[hidelinks,backref=page]{hyperref}
\hypersetup{colorlinks=true,linkcolor=green,citecolor=green,urlcolor=green}

\numberwithin{equation}{section}
\newtheorem{theorem}{Theorem}[section]
\newtheorem{proposition}[theorem]{Proposition}
\newtheorem{lemma}[theorem]{Lemma}
\newtheorem{corollary}[theorem]{Corollary}
\theoremstyle{definition}
\newtheorem{definition}[theorem]{Definition}
\theoremstyle{remark}
\newtheorem{remark}[theorem]{Remark}

\newcommand{\kk}{\Bbbk}
\newcommand{\C}{\mathbb{C}}
\newcommand{\g}{\mathfrak g}
\newcommand{\h}{\mathfrak h}
\newcommand{\m}{\mathfrak m}

\newcommand{\D}{\mathcal D}
\newcommand{\Z}{\mathbb Z}
\newcommand{\N}{\mathbb Z_{\geq0}}
\newcommand{\ad}{\operatorname{ad}}
\newcommand{\End}{\operatorname{End}}
\newcommand{\Cent}{\operatorname{Cent}}
\newcommand{\Wh}{\operatorname{Wh}}
\newcommand{\id}{\operatorname{id}}
\newcommand{\Oe}{\mathcal O_e}

\newcommand{\Uf}{U(\g)_f}

\newcommand{\Proj}{\mathsf P}

\newcommand{\Cchi}{\mathcal C_\chi}
\newcommand{\ba}{\mathbf{a}}
\newcommand{\bb}{\mathbf{b}}

\title[Highest-root localization: a centralizer proof]
{Root vector localization and minimal nilpotent finite \(W\)-algebras}
\author{Genqiang Liu}
\date{\today}
\makeatletter
\@namedef{subjclassname@2020}{%
  \textup{2020} Mathematics Subject Classification}
\makeatother

\begin{document}

\begin{abstract}
Let \(\mathfrak{g}\) be any finite-dimensional complex simple Lie algebra, let $h^\vee$ be its dual Coxeter number, and let \(e\) be
its highest-root vector which is a minimal nilpotent element. In this paper, using a result of Bavula  on
invariants of commuting locally nilpotent derivations, we  show that
\[
 U(\mathfrak{g})_f
 \cong \mathcal{D}^+_{h^\vee-2}\otimes \mathcal{D}_1\otimes U(\mathfrak{g},e),
\]
 where
 the first factor is a  polynomial Weyl algebra of rank $h^\vee-2$, the second factor is a Laurent Weyl algebra of rank one, and  $U(\mathfrak{g},e)$ is the finite $W$-algebra associated with $\mathfrak{g}$ and $e$. As an application, we prove that if $M$ is a simple  infinite-dimensional   $\mathfrak{g}$-module with the minimal GK dimension  on which
 $f-1$ acts locally nilpotently,  then $M\cong P\otimes V$
 for some simple $\mathcal{D}^+_{h^\vee-2}\otimes \mathcal{D}_1$-module $P$ and some finite-dimensional simple  $U(\mathfrak{g},e)$-module $V$. We also give a new proof of the Skryabine quivalence and  Petukhov's birational decomposition of the quotient division ring  of $U(\mathfrak{g})$.
\end{abstract}

\maketitle
\begingroup
\footnotesize
\vspace*{-1.5\baselineskip}
\endgroup

\Needspace{12\baselineskip}
\section{Introduction}

Throughout, \(\C\) is the field of
complex numbers. All associative algebras are unital, and all algebras and
tensor products are over \(\C\).

\subsection{Background and motivation}
Finite \(W\)-algebras are important associative algebras which connect Whittaker modules, Slodowy slices,
primitive ideals of universal enveloping algebras, Yangians, and vertex algebras. A finite $W$-algebra $U(\g, e)$ is an associative algebra associated with a complex simple Lie algebra $\g$ and a nilpotent element $e \in\g$.
From a representation-theoretic perspective, they stem from Kostant's research on Whittaker vectors in a universal Whittaker module attached to a principal nilpotent element \cite{Kostant1978}.  Kostant proved that the corresponding finite \(W\)-algebra is isomorphic to the center \(Z(\mathfrak{g})\) of \(U(\mathfrak{g})\).
Premet constructed the finite \(W\)-algebra attached to an arbitrary
nilpotent element $e$ \cite{Pr1}.
His PBW theorem  identifies its associated Kazhdan graded algebra as
\[
 \operatorname{gr}W(\g,e)\cong\C[\mathcal S_e],
 \qquad \mathcal S_e=e+\g^{f},
\]
where \((e,h,f)\) is an \(\mathfrak{sl}_2\)-triple and
\(\g^{f}=\ker(\ad f)\).
Thus a finite \(W\)-algebra quantizes the Slodowy slice transverse
to a nilpotent orbit.
Skryabin's equivalence reduces generalized Whittaker
\(\g\)-modules to $U(\g,e)$-modules; see \cite[Appendix]{Pr1} and
\cite[Theorem~6.1]{GG}.

In type \(A\), Brundan--Kleshchev realized finite \(W\)-algebras using
truncated shifted Yangians and classified their finite-dimensional simple
modules \cite{BK1,BK2}.
For general type, the connection between finite-dimensional simple representations of
$U(\g, e)$  and primitive ideals of $U(\g)$ was investigated  in \cite{L2,L3,LO,Pr2,Pr4}.
Topley's work described the dimensions of the irreducible components
of the variety of one-dimensional representations in classical type
and related them to the orbit method \cite{T}. Since then, the theory of finite $W$-algebras has been extensively studied by both mathematicians and physicists, see \cite{GG,L10,L1,L2,Pe,PT,RS,VD,We} and references therein.

The theory of \(W\)-algebras also has its origins in mathematical physics.
In two-dimensional conformal field theory, chiral symmetry may extend
the Virasoro algebra by fields of higher conformal spin.
Zamolodchikov's \(W_3\)-algebra is a basic example, with an additional
spin-three field and nonlinear commutation relations \cite{Z}.
A systematic connection with Lie theory is furnished by quantum
Drinfeld--Sokolov reduction, which constructs affine \(W\)-algebras
from affine current algebras using BRST cohomology \cite{FF,KRW}.
At the classical field-theoretic level, Hamiltonian reductions of
Wess--Zumino--Novikov--Witten models lead to generalized Toda theories
with \(W\)-symmetry \cite{FRTW}.
De Boer--Tjin constructed finite \(W\)-algebras by
Poisson reduction and BRST quantization \cite{BT}.
These algebras provide a common language for reduced classical
systems and their quantum symmetries.

The connection between finite and affine \(W\)-algebras is given
by the Zhu algebra construction. If \(W^k(\g,e)\) denotes the universal
affine \(W\)-algebra at level \(k\),
De Sole--Kac proved the associative algebra isomorphism
\(
 \operatorname{Zhu}_{H}\!\left(W^k(\g,e)\right)
 \cong U(\g,e)\),
where \(H\) defines the natural energy grading and
\(\operatorname{Zhu}_{H}\) is the \(H\)-twisted Zhu algebra
\cite[Theorem~5.10]{DeSoleKac}.
The algebra $U(\g,e)$ acts on the lowest-energy spaces of the
corresponding twisted positive-energy modules.
The finite \(W\)-algebra therefore connects the representation theory
of enveloping algebras with that of chiral \(W\)-symmetries.

 When the adjoint orbit of $e$  has the minimal dimension,  the gap between $U(\g, e)$   and $U(\g)$ is particularly small (see \cite{Pr2,Pe}). In this case,  Premet gave generators and relations of $U(\g, e)$ and developed its connection with
Joseph ideals \cite{Pr2}. More recent centralizer constructions of $U(\g,e)$ give tensor product decompositions of $U(\g)$ after
localizing at several root vectors in types \(A\) and \(C\)
\cite{LL,LiuMingjieLi2024}.
In particular finite-dimensional minimal \(W\)-modules were connected with
cuspidal module categories over the corresponding Lie algebras.
Further explicit constructions in type \(A\) are given in
\cite{LiuMingjieLi2025}.
These results motivate a uniform centralizer construction for
localization at a single highest-root vector.

\Needspace{10\baselineskip}
\subsection{The main result}

For $n\geq 1$, let $\D_n^+$ be the rank-$n$ Weyl algebra generated by
$$q_1,\ldots,q_n,p_1,\ldots,p_n$$ with
\[
[p_i,q_j]=\delta_{ij},\qquad
[q_i,q_j]=[p_i,p_j]=0.
\]
For $n=0$, denote  \(\D^+_0=\C\).  Let $\D_1$ be the algebra generated by $x, x^{-1}, D$ such that
$$[D,x]=x,\quad xx^{-1}=x^{-1}x=1.$$

Let \(\g\) be a finite-dimensional complex simple Lie algebra.
Fix a Cartan subalgebra $\h$, and the highest root
$\theta$ of $\g$.  Choose an $\mathfrak{sl}_2$-triple
\[
(e,h,f),\qquad e\in\g_\theta,\quad f\in\g_{-\theta},\quad h=\theta^\vee.
\]
The corresponding $\Z$-grading, defined by the eigenvalues of $\ad h$,
is
\[
\g=\g(-2)\oplus\g(-1)\oplus\g(0)\oplus\g(1)\oplus\g(2),
\qquad \g(2)=\C e,
\]
where
\[
 \g(j)=\{x\in\g:[h,x]=jx\},\qquad
 2s=\dim\g(-1),\qquad d=s+1,
\]
and $ d=\frac12\dim\Oe=h^\vee-1$, \(\Oe\) denotes the adjoint orbit of \(e\).

Our main result of this paper is the following theorem.

\Needspace{13\baselineskip}
\begin{theorem}\label{thm:main}
Let \(U(\g,e)\) be the finite $W$-algebra defined by the generalized Gelfand--Graev
module in Subsection~\ref{subsec:W-algebra}, and $\Uf$
be the localization of $U(\g)$ with respect to the Ore subset
generated by $f$.
Then there is an algebra isomorphism
\begin{equation}\label{eq:main}
 \Uf\cong \D^+_{h^\vee-2}\otimes \D_1\otimes U(\g,e).
\end{equation}
\end{theorem}

The argument is independent of Premet's structural results in \cite{Pr2}.
We use root-space theory, and properties of
Weyl algebras, in particular a result of Bavula \cite{Ba} on
invariants of commuting locally nilpotent derivations.
The finite \(W\)-algebra enters only through its defining
endomorphism algebra. The proof does not use the
PBW theorem of the finite \(W\)-algebra.

We can construct explicit Weyl generators in $\Uf$ and their centralizer.
We first choose root vectors
\[
 z_1,\dots,z_s,w_1,\dots,w_s\in\g(-1),
\]
such that
\begin{equation}\label{eq:Heisenberg}
 [w_i,z_j]=\delta_{ij}f,\qquad
 [z_i,z_j]=[w_i,w_j]=0,\qquad
 [f,z_i]=[f,w_i]=0.
\end{equation}
Choose \(h_0\in\h\) such that
 $\alpha(h_0)\in\Z$ for every root $\alpha\in\Phi$, and
 $\theta(h_0)=1$.
The existence of $h_0$ is proved in Subsection \ref{subsec:cartan-element}.
Suppose  \([h_0,z_i]=r_i z_i\). The explicit Weyl generators are
\begin{equation}\label{eq:explicit-generators}
 q_i=z_i,\quad p_i=f^{-1}w_i,\quad
 D=-h_0+\sum_{i=1}^s r_iq_ip_i,\quad x=f.
\end{equation}
The subalgebras
\begin{equation}\label{eq:B-and-C-intro}
 B=\langle x^{\pm1},D,q_i,p_i:1\leq i\leq s\rangle,
 \qquad C=\Cent_{\Uf}(B)
\end{equation}
satisfy
\[
 B\cong\D^+_{d-1}\otimes \D_1,\qquad \Uf=B\otimes C,\qquad C\cong U(\g,e).
\]
The advantage of this tensor product decomposition is that we can construct explicit generators of $U(\g,e)$ inside the centralizer $C$ using \eqref{eq:retraction}.

\subsection{Relation with some known results}
Define
\[
 \mathcal A_e=
    \bigl(\C[h]\rtimes\langle\Delta\rangle\bigr)\otimes \D^+_s,
 \qquad \Delta h=(h+1)\Delta,
\]
where \(\Delta\) is invertible.
Premet's decomposition \cite[Corollary~5.1(iii)--(iv)]{Pr2}
has the form
\begin{equation}\label{eq:Premet-comparison}
 \Uf\cong(\mathcal A_e\otimes U(\g,e))^{\tau\otimes\sigma}.
\end{equation}
Here \(\tau\) fixes \(h\), sends \(\Delta\) to \(-\Delta\), and
$\tau(p_i)=-p_i, \tau(q_i)=-q_i$.
The map \(\sigma\)  is an involution of $U(\g,e)$
associated with the Dynkin grading.
In this realization, \(f\) maps to \(\Delta^2\), and
\(\mathcal A_\theta\otimes U(\g,e)\) is an extension of the
localization $\Uf$.

Our result gives another realization of the same Ore localization.
The Weyl subalgebra is constructed inside \(\Uf\), and its
centralizer is identified with the defining endomorphism algebra $U(\g,e)=\End_{U(\g)}(Q)^{\mathrm{op}}$.
Our construction does not require a quadratic extension.

Losev's decomposition theorem gives, for a general nilpotent element, an
isomorphism between suitable completions of $U(\g)$ and of a Weyl algebra
tensored with the finite $W$-algebra \cite{L10}.  Our theorem is an
algebraic localization counterpart in the minimal nilpotent setting: it
uses an explicit root vector and produces an isomorphism of
usual Ore localizations, not  of completed topological algebras.
We do not claim that the present theorem is a formal consequence of
Losev's result.

Petukhov studied the Gelfand--Kirillov conjecture for minimal finite
$W$-algebras and the corresponding birational relation with $U(\g)$
\cite{Petukhov}.  Passing from \eqref{eq:main} to quotient division rings
recovers that birational decomposition.  Conversely, a division-ring
isomorphism alone does not imply the Ore-localization isomorphism proved in
this paper.  In this precise sense, the present result refines the
birational picture at the algebra level.

\subsection{Organization}
In Section~\ref{sec:roots}, we fix the root data and prove the Ore property.
In Subsections~\ref{ssec:Weyl} and~\ref{ssec:Laurent}, we establish two Weyl algebra splitting lemmas.
In Subsection~\ref{ssec:construction}, we construct two Weyl algebra factors inside $\Uf$ and their
centralizers.
In Subsection~\ref{ssec:isomorphism},  we show that  the centralizer $C$ is isomorphic to the
minimal finite \(W\)-algebra $U(\g,e)$,  and complete the proof of
Theorem~\ref{thm:main}.
In Section~\ref{sec:applications}, we give some applications of the localization decomposition theorem.

\Needspace{22\baselineskip}
\section{Root data and Ore localization}\label{sec:roots}

In this section, we give the necessary preliminaries including: the good $\Z$-gradings for $e$,
an integral Cartan element, Ore localization, the minimal
finite $W$-algebras.

\subsection{The Heisenberg root spaces}

Throughout this paper, we fix the following notation.
Let $\g=\h\oplus(\oplus_{\alpha\in \Delta}\g_\alpha)$ be a complex simple Lie algebra of rank $r$, with Cartan subalgebra $\h$, root system $\Phi\subset\h^*$, and the set  $\Phi^+$ of  positive roots.
Let $\theta$ be the highest root.  We fix nonzero root vectors $E_\alpha\in\g_\alpha$ and choose $e=E_\theta$, $f\in\g_{-\theta}$ so that
\[
[e,f]=h=\theta^\vee,\qquad [h,e]=2e,\qquad[h,f]=-2f.
\]

Then
\(
\mathfrak s_\theta:=\C e\oplus \C h\oplus \C f 
\cong\mathfrak{sl}_2
\).
We regard $\mathfrak g$ as a finite-dimensional
$\mathfrak s_\theta$-module under the adjoint action.
 By \S 4.3 in \cite{CM}, $e_{\min}:=e$ is a minimal nilpotent element of $\g$. We fix an invariant symmetric bilinear form on $\g$ and normalize it by $(e,f)=1$.

For $j\in \Z$, let $\g(j):=
\{x\in\mathfrak g\mid[h,x]=jx\}$.
By \cite{Pr2}, we have the following $\Z$-grading of $\g$
which is good for $e$.

\begin{lemma}
\label{lem:short-grading}We have the following assertions.
\begin{enumerate}
\item The $\Z$-grading of $\mathfrak g$ defined by
eigenvalues of $\ad h$ is:
\begin{equation}\label{eq:grading}
\mathfrak g
=
\mathfrak g(-2)\oplus
\mathfrak g(-1)\oplus
\mathfrak g(0)\oplus
\mathfrak g(1)\oplus
\mathfrak g(2),
\end{equation}
where
$\mathfrak g(2)=\C e$,
$g(-2)=\mathbb C f$.
\item The form
\[
 \omega(x,y)=(e,[x,y]), \qquad x,y\in\g(-1)
\]
on $\g(-1)$ is nondegenerate and skew-symmetric.
\end{enumerate}
\end{lemma}

\begin{proof}
(1) For a root \(\alpha\), by properties of the root system,
$\frac{2(\alpha,\theta)}{(\theta,\theta)}, \frac{2(\theta,\alpha)}{(\alpha,\alpha)}\in \Z$.
By the
Cauchy--Schwarz inequality, $0\leq \frac{(\alpha,\theta)^2}{(\alpha,\alpha)(\theta,\theta)}\leq 1$. So $\frac{4(\alpha,\theta)^2}{(\alpha,\alpha)(\theta,\theta)}=0,1,2,3,4$.
Since $(\alpha,\alpha)\leq (\theta,\theta)$,
 we must have $\left|\frac{2(\alpha,\theta)}{(\theta,\theta)}\right|\leq 2$. Hence
$|\alpha(h_\theta)|
 =\left|\frac{2(\alpha,\theta)}{(\theta,\theta)}\right|
 \leq2$.
This integer has absolute value \(2\) only when
\(\alpha=\theta\) or \(\alpha=-\theta\), since the root system is reduced
and \(\theta\) has the maximal root length. This proves (1).

(2) By the representation theory of \(\mathfrak{sl}_2\),
\(\ad e:\g(-1)\to\g(1)\) is an isomorphism.
The invariant form pairs \(\g(1)\) and \(\g(-1)\) nondegenerately.
Since
\[
 (e,[x,y])=([e,x],y),
\]
the form \(\omega\) is nondegenerate.

\end{proof}

Define $
\Phi_j=\{\alpha\in\Phi\mid \langle\alpha,\theta^\vee\rangle=j\}$ for $j\in\{-2,-1,0,1,2\}$. Then $\Phi_2=\{\theta\}$, $\Phi_{-2}=\{-\theta\}$.

The following Lemma is crucial for separating the polynomial Weyl algebra $\D^+_{\h^\vee-2}$ from $\Uf$.

\begin{lemma}
\label{lem:short-grading2}
The following assertions hold.
\begin{enumerate}
\item If $\alpha\in\Phi_{-1}$, then $-\theta-\alpha\in\Phi_{-1}$. For  roots  $\alpha, \beta\in \Phi_{-1}$, $[\g_\alpha,\g_\beta]\neq 0$ if and only if $\beta=-\theta-\alpha$.
\item There is a root-vector basis \(z_1,\ldots,z_s,w_1,\ldots,w_s\) of $\g(-1)$
satisfying\begin{equation}
 [w_i,z_j]=\delta_{ij}f,\qquad
 [z_i,z_j]=[w_i,w_j]=0,\qquad
 [f,z_i]=[f,w_i]=0.
\end{equation}
\end{enumerate}
\end{lemma}

\begin{proof}

(1) Since $[\g(-1),\g(-2)]=0$, it follows that $\alpha-\theta$ is not a root.
Then from $\langle\alpha,\theta^\vee\rangle=-1$,  we see that the $\theta$-string through $\alpha$ contains $\alpha+\theta$.  Thus $\alpha+\theta\in\Phi$, equivalently $-\theta-\alpha\in\Phi$.  Its degree is $-1$.

For  roots  $\alpha, \beta\in \Phi_{-1}$, if $[\g_\alpha,\g_\beta]\neq 0$, then $\alpha+\beta$ is a root, and $\alpha+\beta\in \Phi_{-2}$.  Since $\Phi_{-2}=\{-\theta\}$, one has $\alpha+\beta=-\theta$.  Conversely if $[\g_\alpha,\g_{-\theta-\alpha}]=0$, then $[\g_\alpha,\g_\gamma]=0$ and $\omega(\g_\alpha, \g_\gamma)=0$ for any $\gamma\in \Phi_{-1}$. This is a contradiction, since
$\omega(-,-)$ is nondegenerate, by Lemma \ref{lem:short-grading}.

(2) For  roots  $\alpha, \beta\in \Phi_{-1}$, by
 (1), $[\g_\alpha,\g_\beta]\neq 0$ if and only if $\beta=-\theta-\alpha$.
For each pair
\(\{\alpha,-\theta-\alpha\}\), and take \(z_i\in \g_\alpha\) and \(w_i \in\g_{-\theta-\alpha}\)  so that
\([w_i,z_i]=f\).
The other brackets in degree \(-2\) vanish because their roots
cannot sum to \(-\theta\).
Brackets with \(f\) have degree \(-3\), and hence  $[f,z_i]=[f,w_i]=0$.
\end{proof}

\subsection{An integral Cartan element}
\label{subsec:cartan-element}

The following lemma is used  to separate the Laurent Weyl algebra $\D_1$, see the proof of Theorem \ref{thm:centralizer}.

\begin{lemma}\label{lem:integral-Cartan}
There  exists
\(h_0\in\h\) such that
\begin{equation}\label{eq:integral-choice1}
  \theta(h_0)=1, \qquad \alpha(h_0)\in\Z,\quad \forall\,\alpha\in\Phi.
\end{equation}
For the chosen root vectors in Lemma \ref{lem:short-grading2}, we have
\begin{equation}\label{eq:Cartan-weights}
 [h_0,z_i]=r_i z_i,\qquad
 [h_0,w_i]=(-1-r_i)w_i,
\end{equation} where $r_i\in\Z$.
\end{lemma}

\begin{proof}
Choose a simple root \(\beta\) and \(w\) in the Weyl group such that
\(w\beta=\theta\). Let \(\omega_\beta^\vee\in\h\) be the fundamental
coweight, characterized by
\(\gamma(\omega_\beta^\vee)=\delta_{\gamma,\beta}\) for any simple root $\gamma$, and set
\[
 h_0=w\omega_\beta^\vee.
\]
Then \(\theta(h_0)=w\beta(w\omega_\beta^\vee)=\beta(\omega_\beta^\vee)=1\). For any root \(\alpha\),
\(\alpha(h_0)=(w^{-1}\alpha)(\omega_\beta^\vee)\) is an integer.
This proves \eqref{eq:integral-choice1}.
If \(z_i\) has root \(\alpha_i\), then \(w_i\) has root
\(-\theta-\alpha_i\). This proves
\eqref{eq:Cartan-weights}.
\end{proof}

\begin{remark}\label{rem:integrality}
The condition \(\theta(h_0)=1\) gives \([h_0,f]=-f\).
The integer values on all roots will allow powers of \(f\) to remove
every weight in a later centralizer. This is the role of the
integral Cartan choice.
\end{remark}

\subsection{The Ore property and local nilpotence}

By Lemma 4.2 in \cite{Mathieu}, the set $S=\{f^i| i\in \Z_{\geq 0}\}$ is a left and right Ore subset of $U(\g)$.  Let $\Uf$
be the localization of $U(\g)$ with respect to $S$.

\begin{lemma}\label{lem:Ore}
In \(\Uf\), the derivations
\(\ad f,\ad z_i,\ad w_i\) are locally nilpotent.
The algebra \(\Uf\) is a domain.
\end{lemma}

\begin{proof}
The operators \(\ad f,\ad z_i,\ad w_i\) are nilpotent on \(\g\).
This follows directly from the short grading \eqref{eq:grading};
each bracket lowers the degree.
The Leibniz rule extends local nilpotence to \(U(\g)\).

Every element of \(\Uf\) has a right-fraction expression \(uf^{-n}\).
The three derivations \(\ad f,\ad z_i,\ad w_i\)  kill \(f^{-1}\), since \( f, z_i, w_i\)  commute with \(f\). Their local nilpotence therefore
extends to \(\Uf\).
The PBW theorem makes \(U(\g)\) a domain, so these denominators
$\{f^i| i\in \Z_{\geq 0}\}$ are regular (non-zero divisor). Finally, localization at regular elements preserves the domain property.
\end{proof}

\subsection{\texorpdfstring{The minimal finite $W$-algebras}{The minimal finite W-algebras}}
\label{subsec:W-algebra}
 The subspace
\[
 \ell=\operatorname{span}\{w_1,\ldots,w_s\}\subset\g(-1)
\]
is Lagrangian for the form in Lemma~\ref{lem:short-grading}.
The Lie algebra
\[
 \m=\ell\oplus\C f
\]
is abelian.

The following Lemma can found in \cite{Pr2}.

\begin{lemma}\label{lem:dim}$\dim \m=\frac{1}{2}\dim\mathcal{O}_{e}=h^{\vee}-1$
\end{lemma}


The values \(d=h^\vee-1\) are:
\begin{center}
\begin{tabular}{@{}lcc@{}}
\toprule
Type & \(h^\vee\) &  \(d\)\\
\midrule
\(A_\ell\), \(\ell\geq1\) & \(\ell+1\) & \(\ell\)\\
\(B_\ell\), \(\ell\geq2\) & \(2\ell-1\) & \(2\ell-2\)\\
\(C_\ell\), \(\ell\geq2\) & \(\ell+1\) & \(\ell\)\\
\(D_\ell\), \(\ell\geq4\) & \(2\ell-2\) & \(2\ell-3\)\\
\(E_6\) & \(12\) & \(11\)\\
\(E_7\) & \(18\) & \(17\)\\
\(E_8\) & \(30\) & \(29\)\\
\(F_4\) & \(9\) & \(8\)\\
\(G_2\) & \(4\) & \(3\)\\
\bottomrule
\end{tabular}
\end{center}

 Define a linear character  $\chi:\m\rightarrow \C$ by
$\chi(x)=(e,x)$, for $x\in\m$. The character \(\chi\) satisfies
\(\chi(w_i)=0\) and \(\chi(f)=1\). It  defines a one-dimensional $\m$-module $\C_\chi$.

\begin{definition}\cite{Pr1}\label{def:W}
We  define
\[
 Q=U(\g)\otimes_{U(\m)}\C_\chi,\qquad
 U(\g,e)=\End_{U(\g)}(Q)^{\mathrm{op}}.
\]
\end{definition}

Equivalently,
\[
 Q=U(\g)\big/
 (U(\g)(f-1)+\sum_{i=1}^s U(\g)w_i).
\]
The module $Q$ is called the  generalized Gelfand--Graev  module over $\g$.

\section{The tensor product decomposition theorem}
\label{sec:tensor}

In this section, we prove the tensor product decomposition of $\Uf$.

\subsection{Splitting off a polynomial Weyl algebra}
\label{ssec:Weyl}

The following useful result was given by Bavula in \cite{Ba}.
For the convenience of the reader, we provide a proof.

\begin{lemma}\cite[Theorem~2.3]{Ba}\label{lem:Weyl}
Let \(R\) be a unital associative algebra, let \(\delta\) be a locally nilpotent
derivation of \(R\), and let \(v\in R\) satisfy \(\delta(v)=1\).
Then we have
\begin{equation}\label{eq:Taylor-direct-sum}
 R=\bigoplus_{n\geq0}v^n\ker\delta=\C[v]\otimes \ker\delta,
\end{equation}
as vector spaces.
\end{lemma}

\begin{proof}
For $u\in R$, define
\begin{equation}\label{eq:projection}
 \pi(u)=\sum_{n\geq0}\frac{(-v)^n}{n!}\delta^n(u).
\end{equation}
Since $\delta$ is locally nilpotent, this sum is finite.
Then for any $u\in \ker\delta$, we have $\pi(u)=u$. So
\(\pi|_{\ker\delta}=\id\). For $u\in R$, using \(\delta(v^n)=nv^{n-1}\),
\begin{equation}\label{eq:ker}\aligned \delta\pi(u)
&=\sum_{n\geq 1}\frac{-(-v)^{n-1}}{(n-1)!}\delta^n(u)
+\sum_{n\geq0}\frac{(-v)^n}{n!}\delta^{n+1}(u)=0.
\endaligned\end{equation}
Thus \(\pi(R)\subset\ker\delta\).
Since $\sum_{r=0}^{m}\binom{m}{r}(-1)^r=0$ when $m>0$, we have
\begin{equation}\label{eq:inversion}
 \sum_{n\geq0}\frac{v^n}{n!}\pi(\delta^n u)
 =\sum_{n\geq0}\frac{v^n}{n!}
 \sum_{r\geq0}\frac{(-v)^r}{r!}\delta^r(\delta^n u)
 =\sum_{m\geq0} \sum_{r\geq0}^{m}\binom{m}{r}(-1)^r\frac{v^{m}}{(m)!}
\delta^{m}( u)=u,
\end{equation}
it follows that  $R=\sum_{n\geq0}v^n\ker\delta$.

If \(\sum_{n=0}^{N}v^n c_n=0\), with \(c_n\in\ker\delta\),
apply \(\delta^N\) to obtain \(N!c_N=0\).
Characteristic zero gives \(c_N=0\). Downward induction proves
$c_0=\dots=c_{N-1}=0$ and
\eqref{eq:Taylor-direct-sum}.
\end{proof}

By Lemma \ref{lem:Weyl}, we can separate out a Weyl algebra inside an algebra under a natural condition.

\begin{proposition}
\label{lem:finite-separation}
Let $R$ be a unital algebra which contains $\D^+_s$ as a subalgebra,  such that
$\operatorname{ad}p_i$ and
$\operatorname{ad}q_i$
act locally nilpotently on $R$ for every $1\leq i \leq s$.
Then
 we have the following algebra  tensor product decomposition
\begin{equation}\label{eq:Weyl-splitting}R=\D^+_sC_1 \cong \D^+_s\otimes C_1,\end{equation}
where $C_1$ is the centralizer of $\D^+_s$ in $R$.
\end{proposition}

\begin{proof}
We do induction on $s$. The case $s=0$ is immediate. For $s=1$, apply Lemma~\ref{lem:Weyl}
to $\delta=\ad p_1$ and $v=q_1$. This gives
\[
 R=\bigoplus_{a\geq0}q_1^a\ker(\ad p_1).
\]
The derivation $-\ad q_1$ preserves $\ker(\ad p_1)$ because
$[\ad p_1,\ad q_1]=0$, and $(-\ad q_1)(p_1)=1$.
A second application of Lemma~\ref{lem:Weyl} to
the derivation $\ad q_1$ yields
$$\ker(\ad p_1)=\bigoplus_{b\geq0}p_1^b R_1,\qquad
 R_1=\ker_R(\ad q_1)\cap\ker_R(\ad p_1).$$
Thus \begin{equation}\label{eq:k1}
 R=\bigoplus_{a,b\geq0}q_1^ap_1^b R_1.
\end{equation}
For $s>1$, the subalgebra $R_1$ contains $q_i,p_i$ for $2\leq i\leq s$.
The corresponding inner derivations restrict to locally nilpotent
derivations of $R_1$. Their common centralizer inside $R_1$ is precisely
$C_1$. By the induction hypothesis, we have
\begin{equation}\label{eq:m1}
 R_1=\bigoplus_{\alpha,\beta\in\mathbb Z_{\geq0}^{s-1}}
        q_2^{\alpha_2}\cdots q_s^{\alpha_s}
        p_2^{\beta_2}\cdots p_s^{\beta_s}C_1.
\end{equation}
Combining \eqref{eq:k1} and \eqref{eq:m1} proves the direct sum decomposition
\[
 R=\bigoplus_{\ba,\bb\in\mathbb Z_{\geq0}^{s}}
        q^\ba p^\bb C_1,
\]
where $q^{\ba}=q_1^{a_1}\dots q_s^{a_s}, p^{\bb}=p_1^{b_1}\dots p_s^{b_s}$. Thus the multiplication map $\D_s^+\otimes C_1\rightarrow R$ is a vector-space
isomorphism. Since $C_1$ commutes with $\D_s^+$, it is also an algebra
homomorphism, proving \eqref{eq:Weyl-splitting}.
\end{proof}

By Proposition \ref{lem:finite-separation}, every \(u\in R\) has a unique finite expansion
\[
 u=\sum_{\ba,\bb}
 q^{\ba}p^{\bb}
 c_{\ba,\bb},
 \qquad c_{\ba,\bb}\in C_1,
\]
where $q^{\ba}=q_1^{a_1}\dots q_s^{a_s}, p^{\bb}=p_1^{b_1}\dots p_s^{b_s}$.  For later use, we define:
\begin{equation}\label{eq:P-def}
 \begin{aligned}
 \pi_i^p(u)&=\sum_{n\geq0}\frac{(-q_i)^n}{n!}(\ad p_i)^n u,\\
 \pi_i^q(u)&=\sum_{n\geq0}\frac{p_i^n}{n!}(\ad q_i)^n u,\\
 \Proj&=(\pi_s^q\pi_s^p)\cdots(\pi_1^q\pi_1^p).
 \end{aligned}
\end{equation}
By \eqref{eq:ker}, the map $\pi_i^p$ takes values in
$\ker(\ad p_i)$. This kernel is preserved by $\pi_i^q$ because
$[\ad p_i,\ad q_i]=0$. The operators associated with different Weyl
pairs commute. Hence
\begin{equation}\label{eq:P} \Proj(u)\in C_1, \forall\  u\in R.
\end{equation}


\subsection{Splitting off a Laurent Weyl algebra}
\label{ssec:Laurent}

\begin{lemma}\label{lem:Laurent-splitting}
Let \(C_1\) be a unital algebra containing \(x^{\pm1}\) and \(D\) with
\([D,x]=x\). Assume that \(\ad x\) is locally nilpotent and that
the action of \(\ad D\)  on $C_1$ is diagonalizable with integer eigenvalues on \(C_1\).
Then  $C_1=\C[x^{\pm 1}, D]\otimes C\cong \D_1\otimes C$, where
$C=\ker(\ad x)\cap\ker(\ad D)$, and $\C[x^{\pm 1}, D]$ is the subalgebra generated by $x^{\pm 1}, D$.
\end{lemma}

\begin{proof}
Set $T=\ker(\ad x)$. The algebra \(T\) is stable under \(\ad D\).
Indeed, for \(u\in T\), the Jacobi identity gives
\[
 [x,[D,u]]=[-x,u]=0.
\]
The restriction \(\ad D\) to \(T\) gives a weight space decomposition of $T$.
Write \(T=\bigoplus_{m\in\Z}T_m\), with \([D,u]=mu\) for \(u\in T_m\).
The element \(x\) is central in \(T\), is invertible, and $[D,x]=x$.
Hence  \(x^{-m}T_m=T_0=C\).
Thus
\[
 T=\bigoplus_{m\in\Z}x^m C=\C[x^{\pm1}]\otimes C.
\]

Set \(\delta=\ad x\) and \(v=-x^{-1}D\).
Then \(\delta(v)=1\). By Lemma~\ref{lem:Weyl}, we have
\[
 C_1=\bigoplus_{r\geq0}v^rT
  =\bigoplus_{\substack{r\geq0\\m\in\Z}}
     v^r x^m C=\C[x^{\pm 1}, D]\otimes C.
\]
This completes the proof.
\end{proof}

\subsection{The explicit centralizer decomposition of the localization}
\label{ssec:construction}

We now use the root data from Section~\ref{sec:roots}.
In \(\Uf\),  we define \(q_i,p_i,D,x\) by
\begin{equation*}
 q_i=z_i,\qquad p_i=f^{-1}w_i,\qquad
 D=-h_0+\sum_{i=1}^s r_iq_ip_i,\qquad x=f.
\end{equation*}
where
$z_i,w_i$ are defined in Lemma \ref{lem:short-grading2}, and
 $h_0\in \h$ is defined in Lemma \ref{lem:integral-Cartan}.
 
 The following lemma gives Weyl generators of $\D^+_s$ in $\Uf$. 
 
\begin{lemma}\label{lem:fractional-LN}
The relations
\[
 [p_i,q_j]=\delta_{ij},\qquad
 [q_i,q_j]=[p_i,p_j]=0
\]
hold in \(\Uf\). Every \(\ad q_i\) and \(\ad p_i\) is locally nilpotent on $\Uf$.
\end{lemma}

\begin{proof}
These relations follow from \eqref{eq:Heisenberg}, because \(f\)
commutes with all \(z_i,w_i\). By Lemma~\ref{lem:Ore},
the action of \(\ad q_i\) on $\Uf$ is locally nilpotent.

Next we prove the action of $\ad p$ on $\Uf$ is locally nilpotent,  where $p$ denotes \(f^{-1}w_i\).
For \(b\in \Uf\), let \(L_b,R_b\) denote left and right multiplication
by \(b\). Set \(\delta_f=\ad f\) and \(\delta_w=\ad w\).
From
\[
 [f^{-1}w,u]
 =f^{-1}[w,u]-f^{-1}[f,u]f^{-1}w,
\] we have
\begin{equation}\label{eq:fractional-ad}
 \ad(f^{-1}w)
 =
 L_{f^{-1}}\delta_w
 -
 L_{f^{-1}}R_{f^{-1}}R_w\delta_f.
\end{equation}

Since \([f,w]=0\),  $\delta_w$ commutes with $\delta_f$.
All coefficient operators  $L_{f^{-1}}, R_{f^{-1}}, R_w$ in \eqref{eq:fractional-ad} commute with
each other and with both derivations.
Thus, for every \(k\geq0\), the operator binomial theorem gives
\begin{equation}\label{eq:fractional-ad-powers}
 (\ad(f^{-1}w))^k u
 =\sum_{j=0}^{k}(-1)^j\binom kj f^{-k}
   \bigl(\delta_f^j\delta_w^{k-j}u\bigr)f^{-j}w^j.
\end{equation}

For a fixed \(u\), choose \(M,N\geq1\) such that
\(\delta_f^Mu=0\) and \(\delta_w^Nu=0\).
In the binomial expansion of
\((\ad(f^{-1}w))^{M+N-1}u\), every term contains
\(\delta_f^j\delta_w^{M+N-1-j}u\).
Either \(j\geq M\), or \(M+N-1-j\geq N\).
Thus every term is zero. This proves local nilpotence of $\ad p$.
\end{proof}

Let
\begin{equation}\label{eq:C1}
 C_1=\Cent_{\Uf}\langle q_i,p_i:1\leq i\leq s\rangle.
\end{equation}
By Lemma~\ref{lem:Weyl} and Lemma \ref{lem:fractional-LN}, 
\begin{equation}\label{eq:first-splitting}
 \Uf=\C[q_1,\dots,q_s,p_1,\dots,p_s]\otimes C_1\cong \D^+_s\otimes C_1.
\end{equation}
This statement is proved without introducing a finite \(W\)-algebra.

\begin{lemma}\label{lem:D-properties}
The elements \(f^{\pm1}\) and \(D\) belong to \(C_1\), and
\[
 [D,f]=f.
\]
The derivation \(\ad f\) is locally nilpotent on \(C_1\).
The derivation \(\ad D\) is diagonalizable with integer eigenvalues
on \(C_1\).
\end{lemma}

\begin{proof}
By Equation~\eqref{eq:Cartan-weights}, we have
\begin{equation}\label{eq:aqp}
 [h_0,q_i]=r_iq_i,\qquad [h_0,p_i]=-r_ip_i.
\end{equation}
Also,
\[
 [q_ip_i,q_j]=\delta_{ij}q_j,\qquad
 [q_ip_i,p_j]=-\delta_{ij}p_j.
\]
The definition \(D=-h_0+\sum r_iq_ip_i\) therefore gives
\([D,q_i]=[D,p_i]=0\), i.e., $D\in C_1$.
Since \(f\) commutes with $p_i,q_i$,
we have that \(f^{\pm1}\in C_1\), and
\[
 [D,f]=-[h_0,f]=f.
\]
The derivation \(\ad f\) preserves \(C_1\) and is locally nilpotent
by Lemma~\ref{lem:Ore}.

Every \(u\in C_1\) commutes with every \(q_ip_i\). Hence
\begin{equation}\label{eq:D-on-C1}
 [D,u]=-[h_0,u],\quad \forall\ u\in C_1.
\end{equation}
By \eqref{eq:aqp}, the centralizer \(C_1\) is stable under \(\ad h_0\). By \eqref{eq:integral-choice1}, root vectors, Cartan elements, and \(f^{-1}\) are integer eigenvectors
for \(\ad h_0\). So the action of $\ad h_0$ on $\Uf$ is
diagonalizable with integer eigenvalues. The restriction of
$\ad h_0$ on $C_1$ is also
diagonalizable with integer eigenvalues.
Equation~\eqref{eq:D-on-C1} proves the last assertion for $\ad D$.
\end{proof}


Define
\begin{equation}\label{eq:B-and-C}
 B=\langle f^{\pm1},D,q_i,p_i:1\leq i\leq s\rangle,
 \qquad C=\Cent_{\Uf}(B)
\end{equation}

Now we are in a position to prove the main result of this paper. 

\begin{theorem}\label{thm:centralizer}
We have the following  algebra isomorphism
\begin{equation}\label{eq:second-splitting}
 \Uf\cong B\otimes C \cong \D^+_s\otimes \D_1\otimes C.
\end{equation}

\end{theorem}

\begin{proof} By \eqref{eq:first-splitting}, we have \begin{equation}\label{eq:C1-splitting}
 \Uf=\C[q_1,\dots,q_s,p_1,\dots,p_s]\otimes C_1\cong \D^+_s\otimes C_1.
\end{equation}
Apply Lemma~\ref{lem:Laurent-splitting} to \(C_1\), with \(x=f\).
Lemma~\ref{lem:D-properties} verifies all its hypotheses.
Thus
\[
 C_1=\D_1\otimes C.
\]
The algebra \(\D^+_s\) in \eqref{eq:C1-splitting} commutes with
\(C_1\), so combining these two tensor product decompositions proves
\eqref{eq:second-splitting}.
\end{proof}


\subsection{\texorpdfstring{Isomorphism between $C$ and $U(\g,e)$}{Isomorphism between C and U(g,e)}}
\label{ssec:isomorphism}

Recall that $$Q=U(\g)\big/
 (U(\g)(f-1)+\sum_{i=1}^s U(\g)w_i).$$
Let $u_0$ be the image of $1$ in $Q$. For \(u\in U(\g)\), we have
\[
 (f-1)uu_0=[f,u]u_0,
 \qquad
 (f-1)^nuu_0=(\ad f)^n(u)u_0.
\]
Thus \(f-1\) acts locally nilpotently on \(Q\).
Consequently, \(f=1+(f-1)\) acts invertibly on $Q$.
So, $Q$ can be extended uniquely  to be a \(\Uf\)-module.

\begin{lemma}\label{lem:Q-localization}
As a $\Uf$-module,
we have
\begin{equation}\label{eq:Q-localization}
 Q\cong \Uf\big/
 (\Uf(f-1)+\sum_{i=1}^s \Uf p_i ),
\end{equation}
and $\End_{U(\g)}(Q)\cong \End_{\Uf}(Q)$.
\end{lemma}

\begin{proof}
Since $f$ acts invertibly on $Q$, we have $Q\cong \Uf\otimes_{U(\g)} Q$. Let $J=U(\g)(f-1)+\sum_i U(\g)w_i$.
So
$\Uf\otimes_{U(\g)}(U(\g)/J)\cong \Uf/\Uf J$,
and
\(Q\cong \Uf/(\Uf(f-1)+\sum_i \Uf w_i)\).
Since \(w_i=fp_i\) and \(f\) is invertible,
\(\Uf w_i=\Uf p_i\). This proves \eqref{eq:Q-localization}.

An endomorphism commuting with \(f\) also commutes with its inverse.
Therefore every \(U(\g)\)-endomorphism is a \(\Uf\)-endomorphism.
\end{proof}

Define the following $B$-module
\[
 B_1=B\big/(B(f-1)+\sum_{i=1}^s Bp_i).
\]
All generators of the defining left ideal belong to \(B\).
By the factorization \(\Uf=B\otimes C\),
\begin{equation}\label{eq:Q-factor}
 Q\cong B_1\otimes C.
\end{equation}

Let \(v_0\) be the image of \(1\) in $B_1$. For any \(\g\)-module \(M\) on which $f-1$ and $w_i$ act locally nilpotently for all $i$, we define the Whittaker space
\[
 \Wh(M)=\{v\in M| (f-1)v=w_iv=0\text{ for all }i\}.
\]

\begin{lemma}\label{lem:polynomial-model}
In \(Q=B_1\otimes C\),
we have
 $\Wh(Q)=v_0\otimes C$.
\end{lemma}

\begin{proof}

By the construction of $B_1$, there is a \(B\)-module isomorphism
\[
 B_1\cong\C[y,q_1,\ldots,q_s]
\]
where $y= \partial=f^{-1}D$. On $B_1$, the actions of
$f, p_i$ are given by
\begin{equation}\label{eq:polynomial-action}
 f=1-\partial_y,\qquad p_i=\partial_{q_i}.
\end{equation}

On \(B_1\otimes C\), the equation \((f-1)v=0\) is
\(\partial_yv=0\).
Since \(w_i=fp_i\) and \(f\) is invertible, \(w_iv=0\) is equivalent to
\(\partial_{q_i}v=0\).
So the joint kernel of these derivatives is \( v_0\otimes C\).
This completes the proof.
\end{proof}

\begin{theorem}
\label{thm:W-identification}
$C\cong\End_{U(\g)}(Q)^{\mathrm{op}}
\cong U(\g,e).$
\end{theorem}

\begin{proof}
By Lemma \ref{lem:polynomial-model}, the joint kernel of $f-1$ $w_i$ on
$Q=B_1\otimes C$ is $v_0\otimes C$.  Therefore a
$(B\otimes C)$-endomorphism of $Q$ is determined by
\[
v_0\otimes1\mapsto v_0\otimes c,
\] where $c\in C$.
  Thus
\[
\End_{U(\g)}(Q)\cong \End_{B\otimes  C}(B_1\otimes C)\cong \End_{C}(C)
\cong C^{\mathrm{op}}.
\]
Therefore
$C\cong\End_{U(\g)}(Q)^{\mathrm{op}}
\cong U(\g,e)$.
\end{proof}

Combining Theorems \ref{thm:centralizer} and \ref{thm:W-identification}, we obtain the main result in this paper.

\begin{theorem}\label{thm:main1}
Let $\Uf$
be the localization of $U(\g)$ with respect to the Ore subset
generated by $f$.
Then there is an algebra isomorphism
\begin{equation}\label{eq:main1}
 \Uf\cong \D^+_{h^\vee-2}\otimes \D_1\otimes U(\g,e).
\end{equation}
\end{theorem}

\begin{proof}
Theorems~\ref{thm:centralizer} and~\ref{thm:W-identification} give
\(\Uf\cong\D^+_s\otimes \D_1\otimes U(\g,e)\).
Since  $\dim \m=\frac{1}{2}\dim\mathcal{O}_{\mathrm{min}}=h^{\vee}-1$,
we have $s+1=h^{\vee}-1$. This proves \eqref{eq:main1}.
\end{proof}

\begin{remark} We  conjecture that this type of tensor product decomposition should be valid for a broader class of Lie algebras and Lie superalgebras. For Cartan type Lie algebras, we can see \cite{LZ24}.
\end{remark}

\section{Applications}\label{sec:applications}

By Theorem \ref{thm:main1}, there is an algebra isomorphism
\begin{equation}\label{eq:main22}
 \Uf\cong \D^+_{h^\vee-2}\otimes \D_1\otimes U(\g,e).
\end{equation}
In this section, we give several applications of this decomposition theorem.

\subsection{\texorpdfstring{Representations of $\g$ with the minimal GK dimension}{Representations with the minimal GK dimension}}

In this subsection, we will use simple holonomic $\D^+_{h^\vee-2}$-modules and finite-dimensional simple $U(\g,e)$-modules to study infinite-dimensional simple $\g$-modules with the smallest GK dimension.

\begin{lemma} Let $M$ be an infinite-dimensional simple $U(\mathfrak g)$-module. Then $\operatorname{GKdim}M
\geq h^{\vee}-1$.
\end{lemma}

\begin{proof}
Set
$I=\operatorname{Ann}_{U(\mathfrak g)}M$.
By the generalized Bernstein inequality \cite[Section~1.2]{Losev2017},
\[
2\operatorname{GKdim}M
\geq
\operatorname{GKdim}\bigl(U(\mathfrak g)/I\bigr)
=
\dim\mathcal V(I).
\]
Since $I$ is primitive, its associated variety $\mathcal V(I)$ \cite{Joseph1,Joseph2} is the closure of a
nilpotent orbit:
$\mathcal V(I)=\overline{\mathcal O}$.
The module $M$ is infinite-dimensional, so $\mathcal O\neq\{0\}$.
Consequently,
$\dim\mathcal O\geq\dim\mathcal O_{\min}$.
By \cite{W1},
$\dim\mathcal O_{\min}=2h^\vee-2$.
Therefore
$2\operatorname{GKdim}M
\geq
2h^\vee-2$,
which gives
$\operatorname{GKdim}M\geq h^\vee-1$.
\end{proof}

\begin{definition}\label{def:minimal-representation}

We call an infinite-dimensional simple $U(\mathfrak g)$-module $M$
 a \emph{minimal representation} of
$\mathfrak g$ if
\[
\operatorname{GKdim}M
=
\frac{1}{2}\dim\mathcal O_{\min}=h^{\vee}-1.
\]
\end{definition}

The following tensor product lemma is standard; see
\cite[Lemma~3.6]{XL}.

\begin{lemma}\label{tensor}
Let $A$ and $B$ be two unital associative algebras, with \(A\)
countable-dimensional.
\begin{enumerate}
\item If $M$ is a simple module over $A\otimes B$ that contains a simple $A=\C\otimes A$ submodule $V$, then $M \cong V\otimes W$ for a simple $B$-module $W$.
\item If $V$ and $W$ are simple modules over $A$ and $B$ respectively, then $V\otimes W$ is a simple module over $A\otimes B$.
\end{enumerate}
\end{lemma}

Next using the isomorphism \eqref{eq:main22},
we classify a family of simple minimal representations of
$\mathfrak g$.

\begin{theorem}Suppose that $M$ is a  simple minimal $U(\mathfrak g)$-representation on which $f-1$ acts locally nilpotently. Then there exist a simple $\D^+_{h^\vee-2}$-module $V$ of $\operatorname{GKdim} V=h^\vee-2$, and a finite-dimensional simple
$U(\g,e)$-module $N$ such that $M\cong V\otimes \D/\D (x-1)\otimes N$.
\end{theorem}

\begin{proof} By the isomorphism $\Uf\cong \D^+_{h^\vee-2}\otimes \D_1\otimes U(\g,e)$, $M$ can be viewed as a simple module over
$\D^+_{h^\vee-2}\otimes \D_1\otimes U(\g,e)$   on which $x-1$ acts locally nilpotently. Choose nonzero $v\in M$ such that $xv=v$. We can check that $\D_1 v=\D_1/\D_1(x-1)$ is a simple $\D_1$-module whose GK dimension is $1$. By Lemma \ref{tensor}, there is a simple $\D^+_{h^\vee-2}\otimes U(\g,e)$-module $K$ such that $M\cong  \D_1/\D_1(x-1)\otimes K$. The GK dimension of $K$ is $h^\vee-2$.  So
$K$ is a holonomic module over the Weyl algebra  $\D^+_{h^\vee}$. Since a holonomic module over  $\D^+_{h^\vee-2}$ has a finite composition length, restricted to $\D^+_{h^\vee}$,
$K$ has a simple $\D^+_{h^\vee-2}$-submodule $V$. Again  by Lemma \ref{tensor}, there is a simple $W(\g,e)$-module $N$ such that $K\cong V\otimes N$. Since
$\operatorname{GKdim}K=\operatorname{GKdim} V=h^\vee-2$,  we have $\operatorname{GKdim} N=0$. Hence $N$ is finite-dimensional.
\end{proof}

\begin{theorem}Suppose that $\g$ is of type $A$, $D$ or $E$, and $M$ is a  simple minimal $U(\mathfrak g)$-representation. Then
there is an automorphism $\sigma$ of $\g$ such that
$\Uf\otimes_{U(\g)} M^{\sigma}\cong K\otimes N$, for some simple $\D^+_{h^\vee-2}\otimes \D_1$-module $K$ and some finite-dimensional simple $U(\g,e)$-module $N$.
\end{theorem}

\begin{proof} By \cite[Lemma~2.4]{LZ}, there exists a root vector $e_\alpha$ acting injectively on $M$. Since $\g$ is of type $A$, $D$ or $E$, there is an automorphism $\sigma$ of $\g$ such that $\sigma(f)=e_{\alpha}$. So $f$ acts injectively on $M^\sigma$ which is obtained by twisting $M$ by $\sigma$. Then $\Uf\otimes_{U(\g)} M^{\sigma}$ is a simple
$\Uf$-module of GK dimension $h^\vee-1$. Consequently this theorem follows from Lemma \ref{tensor}.
\end{proof}

It should be mentioned that finite-dimensional simple
$U(\g,e)$-modules were described in \cite{Pr2}.

\subsection{A new proof of the Skryabine quivalence}

In this subsection, we identify $U(\g,e)$ with the centralizer $C$ of
$B$ in $\Uf$.

\begin{lemma}\label{lem:module-Taylor}
If \(M\) is a $D^+_1$-module on which  \(p\) is locally nilpotent.
Then
\[
 \pi_M(v)=\sum_{n\geq0}\frac{(-q)^n}{n!}p^nv
\]
takes values in \(\ker p\), is the identity there, and gives
\[
 v=\sum_{n\geq0}\frac{q^n}{n!}\pi_M(p^nv),\qquad
 M=\bigoplus_{n\geq0}q^n\ker p.
\]
\end{lemma}

\begin{proof}The proof is similar as that of Lemma \ref{lem:Weyl}.
\end{proof}

\begin{definition}\label{def:Wh-category}
Let \(\Cchi\) be the category of all \(U(\g)\)-modules on which
\(x-\chi(x)\) acts locally nilpotently for every \(x\in\m\).
\end{definition}

For every object $M\in \Cchi$,
\(f-1,w_1,\ldots,w_s\) act on it locally nilpotently.
Then $M$ can be equipped with a \(\Uf\)-module structure, since  \(f\) acts bijectively on $M$.
Since \(p_i=f^{-1}w_i\) and \([f,w_i]=0\), local nilpotence of
\(w_i\) is equivalent to local nilpotence of \(p_i\).
Since $[U(\g,e),f]=[U(\g,e),w_i]=0$,  \(U(\g,e)\) preserves \(\Wh(M)\).

\begin{theorem}\cite{S}\label{thm:elementary-Skryabin}
The functors
\begin{equation}\label{eq:Wh-equivalence}
 \begin{aligned}
 U(\g,e)\text{-}\mathrm{Mod}&\longrightarrow\Cchi,
       &N&\longmapsto Q\otimes_{U(\g,e)} N,\\
 \Cchi&\longrightarrow U(\g,e)\text{-}\mathrm{Mod},
       &M&\longmapsto\Wh(M)
 \end{aligned}
\end{equation}
are mutually inverse equivalences. Explicitly, with
\(B_1=\kk[y,q_1,\ldots,q_s]\) as in
Lemma~\ref{lem:polynomial-model},
\begin{equation}\label{eq:Wh-module-factor}
 Q\otimes_{U(\g,e)} N\cong B_1\otimes N,\qquad
 M\cong B_1\otimes\Wh(M).
\end{equation}
\end{theorem}

\begin{proof}
On \(M\in\Cchi\), use the commuting pairs
\[
 p_0=f-1,\quad q_0=-f^{-1}D,\qquad
 p_i=f^{-1}w_i,\quad q_i=z_i, \quad 1\leq i\leq s.
\]
They satisfy \([p_i,q_j]=\delta_{ij}\), and different pairs commute.
All \(p_i\) are locally nilpotent.
Apply Lemma~\ref{lem:module-Taylor} repeatedly, we have
\begin{equation}\label{eq:module}
 M=\bigoplus_{r\geq0,\,\ba\in\N^s}
       q_0^r q^{\ba}\Wh(M).
\end{equation}

Define
\[
 B_1\otimes\Wh(M)\rightarrow M,
 \qquad q_0^r q^{\ba}\otimes v
 \mapsto q_0^r q^{\ba}v.
\]
The expansion \eqref{eq:module} proves bijectivity. The Weyl commutators show that
\(p_i\) acts as \(\partial_{q_i}\), and
\(f\) acts as \(1-\partial_y\).
Thus the map is a \(B\)-module homomorphism. It is a \(U(\g,e)\)-homorphism because \(U(\g,e)\)
commutes with all the displayed Weyl generators.
It is therefore a \(U(\g)\)-homomorphism.

Conversely, \(B_1\otimes N\) belongs to \(\Cchi\), and its joint
kernel of \(f-1,p_1,\ldots,p_s\) is \(1\otimes N\).
Since \(Q=B_1\otimes U(\g,e)\) with its right \(U(\g,e)\)-action,
\(Q\otimes_{U(\g,e)}N=B_1\otimes N\). These are the two natural inverse maps.
\end{proof}

\subsection{Recovery of Petukhov's birational decomposition}

\begin{corollary}\cite{Petukhov}\label{cor:quotient-field}
Let $Q(-)$ denote the quotient division ring of a domain.  Then
\[
Q(U(\g))\cong
Q\!\left(\D^+_{h^\vee-1}\otimes U(\g,e_{\min})\right).
\]
\end{corollary}

\begin{proof}
The  associated graded algebra of $U(\g,e_{\min})$ with respect to the Kazhdan filtration is the
coordinate ring of the Slodowy slice \cite{Pr1}, hence it is a domain; consequently
$U(\g,e_{\min})$ is a Noetherian domain.  Passing to
$U(\g)_f$ does not change the quotient division ring.  By
\eqref{eq:main22}, its quotient division ring is that of
\[
\D^+_{h^\vee-2}\otimes\D_{1}\otimes U(\g,e_{\min}).
\]
The first two Weyl factors form a localization of
$\D^+_{h^\vee-1}$ and hence have the same quotient division ring as
$\D^+_{h^\vee-1}$.  The assertion follows.
\end{proof}

Thus the algebraic localization theorem recovers, after passage to
quotient division rings, the birational relation studied by Petukhov
\cite{Petukhov}.
Therefore the birational structure of $U(\mathfrak g)$ splits into a
standard Weyl part and a residual part determined by the minimal finite
$W$-algebra.  In particular, it is closely related to
Gelfand--Kirillov-type problems for both algebras.
We can show that if an analogue of the Gelfand-Kirillov conjecture
holds for such $W(\g,e_{\min})$, then it holds for the universal enveloping algebra $U(\g)$.
The solution of the  Gelfand-Kirillov conjecture for Lie algebras of type $A$  was
settled by Gelfand and Kirillov  in \cite{GK1}.
Premet showed in \cite{Pr3} that the Gelfand--Kirillov conjecture fails for $U(\g)$ if $\g$ is simple and $\g$ is
not of type $A_n, C_n$ or $G_2$.

\vspace{4mm}

 \noindent G.L.: School of Mathematics and Statistics,
Henan University, Kaifeng 475004, China. Email: liugenqiang@henu.edu.cn

\vspace{0.2cm}

\end{document}